\documentclass[12pt]{article}
\usepackage{psfig}
\begin{document}
 
\bibliographystyle{plain}

\author{Garrison W. Greenwood\thanks{email:
greenwd@ee.pdx.edu}} 

\date{Dept. of Elec.~\&
Computer Engineering, Portland State University, Portland, OR
97207}

\title{On the Isomorphic Embedding of Rectangular 
Grids in \emph{n}-cubes
}

\maketitle

\baselineskip 21pt

\begin{abstract}
All previously published work on isomorphic 
grid embeddings into $n$-cubes
has been restricted to binary $n$-cubes.
This paper describes a straightforward method for embedding a 
$A\times B$ grid isomorphically into a $k$-ary $n$-cube with $k>2$.  
\end{abstract}

Suppose you have meticulously assigned the tasks of an algorithm 
to these processors taking particular care to ensure the task
assignments will require information to travel minimal distances
over the network.  Now, now you are asked to run this
algorithm on a different parallel processing system that has a
different network architecture.  How can you achieve the
previous performance with minimal reprogramming effort?  

The answer depends on whether
or not an \emph{isomorphic embedding} exists between the two
network architectures.
Consider a graph $G=(V,E)$ with a set of vertices or nodes $V$
and a set of nondirectional edges $E$ which connect pairs of vertices.  
With respect to the problem above, the parallel processing system
can be depicted as a graph where each processor is a vertex
and the edges are interconnection links that form the
communications network.  
The embedding of one graph
(called the target graph) into another graph (called
the host graph) assigns each vertex in the target graph to a
vertex in the host graph.  Similarly, edges in the target
graph will overlay edges in the host graph.  

Isomorphic embeddings do not always exist between a target and
host graph, and even if they do, finding them is NP-hard.  
Nevertheless, isomorphic embeddings are of particular interest to
the computer science field because of their strong implications 
for parallel processing:  given two parallel machines 
$M_1$ and $M_2$ with network topologies $T_1$ and $T_2$ respectively, 
if $T_2\mapsto T_1$, then any algorithm that runs in $N$ steps 
on $M_1$ will likewise run in $N$ steps on 
$M_2$\footnote{ ``$\mapsto$'' denotes an isomorphic embedding.}.
In other words, if an isomorphic embedding exists between the 
graphs depicting the two parallel processing systems, then little 
or no reprogramming effort is required.  In fact, the previously
attained performance on the old system will be duplicated in the new 
system.

\begin{figure}[htbp]
\centerline{\psfig{figure=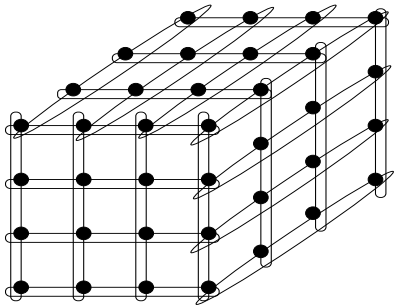,height=2.0in,width=2.0in}}
\caption{A $k$-ary $n$-cube with $k=4$ and $n=3$.  Hidden nodes and 
edges are not shown to preserve clarity.}
\label{ncube}
\end{figure}

One particularly versatile network topology is the $n$-cube.
$k$-ary $n$-cubes are graphs with $n$ dimensions and $k$ 
nodes in each direction.  
(An example of a 4-ary 3-cube is shown in 
Figure \ref{ncube}.)  
It is well known that binary hypercubes (i.e., 2-ary $n$-cubes) 
are efficient architectures for executing parallel 
algorithms.  Their reasonable tradeoff between number of processors and 
interconnectivity makes them ideally suited for solving linear algebra 
and graph-theoretical problems \cite{dekel81}.
The performance of the more general \emph{k}-ary \emph{n}-cube has
been analyzed by Dally \cite{dally90b}.

All known previous work involving rectangular grid embeddings 
into $n$-cubes has been restricted to 2-ary $n$-cube host
graphs.  In certain restricted cases it is possible to 
isomorphically embed a rectangular grid into a 2-ary $n$-cube 
\cite{saad88}---universally applicable techniques require dilation 
costs greater than unity \cite{bett92,ho87}. 

This paper describes a straightforward method for embedding a 
$A\times B$ grid 
isomorphically into a $k$-ary $n$-cube with $k>2$.  
The isomorphic 
embedding described shortly can always be 
accomplished with the mild restriction of $A\le k$.  Neither $A$ 
nor $B$ are required to be an integer power of 2. 
Moreover, the technique does not require embedding into 
an intermediate graph (e.g., a binary hypercube or butterfly
network).

Let $b_{1}b_{2}\ldots b_n$ be a $n$-bit binary string 
where $b_i\in\{0,1\}$.  There are 2$^n$ unique binary patterns 
that can be formed with $n$-bit binary strings.  A sequence of 
length $L$ contains $L$ $n$-bit binary strings no two of which 
are identical.  

Let $A$ and $B$ be integers.  The objective is isomorphically 
embed a $A\times B$ grid into a $k$-ary $n$-cube, where $A\le k$ 
and $B\le k^{(n-1)}$.  For the moment, assume $A=k$ and $k$ is 
an integer power of two---restrictions that will be shortly 
removed.  Each node in the $k$-ary $n$-cube is labeled with 
a $n\log_2 k$ bit binary label.  The labeling is done in a Gray code manner such 
that any two nodes connected by an edge differ in only one bit 
position.  

The binary label associated with each node in the $k$-ary $n$-cube 
can be partitioned into two parts 
\[
\{\underbrace{ b_1 \ldots b_{r}}_{\lceil \log_2 A \rceil} \  \underbrace{b_{r+1} 
\ldots b_{n\log_2 k}}_{\lceil\log_2 B\rceil} \}
\]
Thus any point $(x,y)$ in the 2-dimensional grid is given by using 
the $\lceil\log_2 A\rceil$ most significant bits to define the $x$ 
coordinate, and the $\lceil\log_2 B\rceil$ least significant bits to 
define the $y$ coordinate.  Any two adjacent nodes in the cube will 
be adjacent in the grid---a  property enforced by the Gray code 
labeling.  

\begin{figure}[htbp]
\centerline{\psfig{figure=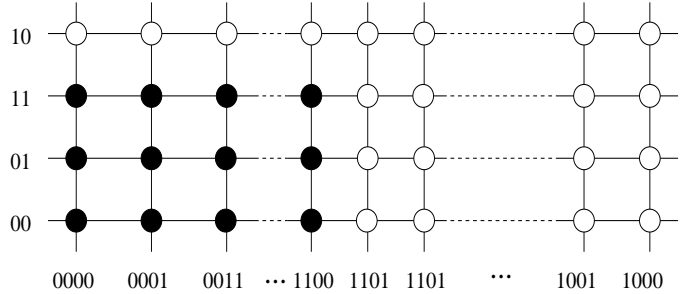,height=1.5in,width=3.5in}}
\caption{An $A\times B$ grid embedding into a 4-ary 3-cube with $A=3$ and $B=9$.  
The placement of the rows of the grid are determined by the 
$\lceil\log_2 A\rceil$ most significant bits and the columns
are labeled with the remaining bits of the binary label.  The darkened
nodes are used in the grid embedding.}
\label{latt}
\end{figure}

In practice, it is not necessary to have $A$ and $B$ as integer 
powers of 2; the use of integer ceilings in the label partitioning 
insures there are a sufficient number of bits.  For example, consider 
a case where $A=3$, $B=9$, and the grid is to be embedded into a 
4-ary 3-cube.  Note that $A\le 4$ and $B\le 4^2=16$, which means the 
cube is large enough to contain the grid.  
Each node in the $n$-cube has a 
3$\;\log_2 4= 6$-bit binary label.  Then $\lceil\log_2 A\rceil=2$ bit 
positions in each binary label are reserved for grid row assignments, 
and $\lceil\log_2 B\rceil=4$ bit positions are used for grid column 
assignments. Incrementing the respective sets of bits in
a Gray code manner will identify the ultimate assignments.
Figure \ref{latt} shows a $3\times 9$ grid embedding in a 4-ary 
3-cube.

In effect, this grid embedding technique ``unrolls'' a high 
dimensional cube into a 2-dimensional grid.  This unrolling does 
break edges in the $n$-cube and so some neighbor relationships are 
lost.  Nevertheless, the grid embedding is isomorphic.  If $k$ is 
not an integer power of 2, set $\tilde k = 2^{\lceil \log_2 k \rceil}$, 
and then isomorphically embed the grid into a $\tilde k$-ary $n$-cube.


\begin{thebibliography}{99}
\bibitem{bett92} S. Bettayeb and Z. Miller, Embedding grids into
hypercubes, \emph{J. Comp. \& Sci.} 45 (1992) 340-266.
\bibitem{dally90b} W. Dally, Performance analysis of k-ary n-cube
interconnection networks, \emph{IEEE Trans. Comp.} 39 (1990) 775-785.
\bibitem{dekel81} E. Dekel, D. Nassimi, and S. Sahni, Parallel matrix
and graph algorithms, \emph{SIAM J. Comp}. 10 (1981) 657-675.
\bibitem{ho87} C. Ho and S. Johnsson, On the embedding of arbitrary
meshes in boolean cubes with expansion two dilation two, \emph{Proc.
Int'l Conf. on Para. Proc.} (1987) 188-191.
\bibitem{saad88} Y. Saad and M. Schultz, Topological properties of
hypercubes, \emph{IEEE Trans. Comp.} 37 (1988) 867-871.
\end{thebibliography}
\end{document}